\documentclass{article}
\usepackage{amsmath}
\usepackage{amssymb}
\usepackage{amsthm}
\theoremstyle{plain}
\newcommand{\field}[1]{\mathbb{#1}}

\newcommand{\R}{\field{R}}

\newcommand{\beeq}{\begin{equation}}
\newcommand{\eneq}{\end{equation}}
\newcommand{\la}{\label}
\newcommand{\pref}[1]{(\ref{#1})}
\newcommand{\ba}{\begin{array}}
\newcommand{\ea}{\end{array}}

\newcommand{\be}{\begin{equation}}
\newcommand{\ee}{\end{equation}}
\newcommand{\bea}{\begin{eqnarray}}
\newcommand{\eea}{\end{eqnarray}}
\newcommand{\nn}{\nonumber}

\newcommand{\intmean}{{\int\hspace{-10.6pt}-}}
\newcommand{\smallintmean}{{\int\hspace{-8.6pt}-}}
\newcounter{const}

\newtheorem{thm}{Theorem}[section]

\newtheorem{lem}[thm]{Lemma}

\newtheorem{rem}[thm]{Remark}

\title{
$C^{1,\alpha}$-Regularity of energy minimizing maps
from a 2-dimentional domain into a Finsler space
}
\author{
Atsushi Tachikawa
\thanks{
This research was partially supported by the Ministry of Education, Science, Sports and Culture, Grant-in-Aid for Scientific Research (C), 22540207}
\\
\small Department of Mathematics, Faculty of Science and Technology,\\
\small Tokyo University of Science, Noda, Chiba, 278-8510, Japan,\\
\small e-mail:tachikawa$\_$ atsushi@ma.noda.tus.ac.jp
}

\date{}

\begin{document}
\maketitle
\noindent
2010\textit{MSC}:
Primary 49N60, 58E20; Secondary
35B65, 53C60.
\begin{abstract}
We show $C^{1,\alpha}$-regularity for energy minimizing maps from
a 2-dimensional Riemannian manifold into a Finsler space $(\R^n, F)$
with a Finsler structure $F(u,X)$.
\end{abstract}
\section{Introduction}
Let $N$ be an $n$-dimensional $C^\infty$-manifold and 
$TN$  its tangent bundle.
We write each point in $TN$ as $(u,X)$ with $u\in N$ and
$X\in T_u N$. 
We put 
 \[
 TN\setminus 0 := \{(u,X)\in TN~;~ X\neq 0\}.
 \]
$TN\setminus 0$ is called the \it slit tangent bundle \rm of $N$.
A \em Finsler structure \rm of $N$
is a function $F \colon TN \to [0,\infty)$
with the following properties:
\begin{description}
  \item {(F-1) \bf{Regularity:}} $F \in C^\infty(TN\setminus 0)$.
  \item {(F-2) \bf{Positive homogeneity:}} 
   $F(u,\lambda X) = \lambda F(u,X)~~\text{for all}~\lambda\geq0.$
  \item {(F-3) \bf{Convexity:}} The Hessian matrix of $F^2$ with respect to $X$ 
\[ 
     (f_{ij} (u,X)) = \left( \frac{1}{2}
     \frac{\partial^2 F^2 (u,X)}{\partial X^i \partial X^j}
     \right)
\]
is positive definite at every point $(u,X) \in TN\setminus 0$.
\end{description}
We call the pair $(N,F)$ a \em Finsler manifold, \rm
and $(f_{ij})$ the \it fundamental tensor \rm of $(N,F)$.
Since $F$ is positively homogeneous of degree 1, we can see that
the coefficients of the fundamental tensor are positively homogeneous
of degree 0;
\begin{equation}\la{hom-h}
f_{ij}(u,\lambda X)=f_{ij}(u,X), ~~~\lambda >0.
\end{equation}
Moreover, since $F^2$ is homogeneous of degree 2, 
using Euler's theorem for homogeneous functions,
we have
\begin{equation}\la{Euler}
F^2(u,X)=f_{ij}(u,X)X^iX^j.
\end{equation}

For maps  between Finsler manifolds 
P. Centore \cite{cen00} defined the \it energy density \rm  
by using of the integral mean
on the indicatrix of each point on the source manifold. 
According to his definition
we define the \it energy density \rm $e_C(u)$ of a map $u$ 
from a Riemannian into a Finsler manifold as follows.
Let $(M,g)$  be a smooth Riemannian $m$-manifold and  $(N,F)$
a Finsler $n$-manifold. Let $I_xM$ be the indicatrix of $g$ at $x \in M$, namely,
\[
I_xM := \{\xi \in T_xM ; \|\xi\|_g \leq 1\}.
\]
For a $C^1$-map $u:M \to N$
and a domain $\Omega \subset M$,
we define the \it energy density \rm
$e_C(u)(x)$ of $u$ at $x\in M$ and the \it energy \rm on $\Omega$  
$E_C(u;\Omega)$ by
\begin{align}
    & e_C(u)(x) := \intmean_{I_xM}(u^\ast F)^2(\xi)d\xi=
    \frac{1}{\int_{I_xM}d\xi}\int_{I_xM} (u^\ast F)^2(\xi)d\xi
   \la{e-dens}\\
   & E_C(u;\Omega):= \int_\Omega e_C(u)(x) d\mu . \la{e-ful}
\end{align}
Here and in the sequel, $\smallintmean$ denotes the integral mean,
$u^\ast F$ the pull-back of $F$ by $u$, 
and $d\mu$ the measure deduced from $g$.
We call (weak) solutions of the
Euler-Lagrange equation of the energy  \it (wakly) harmonic
maps. \rm

Concerning harmonic maps from a Finsler manifold into a Riemannian manifold,
see, for example, H. von der Mosel and S. Winklmann \cite{vonwin09}.

Let us  take an orthonormal frame 
$\{e_\alpha\}$ for the tangent bundle $TM$ of $M$, given in local coordinates by
\[
e_\alpha = \eta^\kappa_\alpha(x) \frac{\partial}{\partial x^\kappa}, ~~1\leq \alpha \leq m.
\]
Using $\{e_\alpha\}$, we identify each $I_xM$
at $x \in M$ with the unit Euclidean $m$-ball $B^m$.
Then, by virtue of the identity
\[
 g^{\kappa\nu}(x)= \eta^\kappa_\alpha(x)\delta^{\alpha\beta}\eta^\nu_\beta(x),
\] 
we can write $E_C$ as
\begin{align}
	& E_C (u;\Omega) \nn\\
	= & \int_\Omega \left( \frac{1}{|B^{m}|} 
	\int_{B^{m}} f_{ij}(u(x),du_x(\xi) ) \xi^\kappa \xi^\nu d\xi \right) 
	\eta^\alpha_\kappa \eta^\beta_\nu 
	D_\alpha u^i D_\beta u^j \sqrt{g}dx,  \la{C-energy}
\end{align}
where $D_\alpha u^i= \partial u^i/\partial x^\alpha$
and $g= \det(g_{\alpha\beta})$.
(cf.  \cite{nis08}.)
Although the terms in parentheses are not defined at points $x$ where $du_x=0$,
we can define them to be arbitrary numbers without changing the values of the integrands
$
(.....)\eta^\alpha_\kappa \eta^\beta_\nu 
	D_\alpha u^i D_\beta u^j,
$
because the integrands are equal to $0$, being independent
on the values of $f_{ij}$ when $du_x=0$.  
So, here and in the sequel, we regard $f_{ij} (u,X)$  as being
defined also for $X=0$. 

As in \cite{tac09}, let us put
\begin{align}
	& E^{\alpha\beta}_{ij}(x,u,p)\nn\\
	= & \left( \frac{1}{|B^{m}|} 
	\int_{B^{m}} f_{ij}(u(x),p\xi) ) \xi^\kappa \xi^\nu d\xi \right) 
	\eta^\alpha_\kappa (x) \eta^\beta_\nu (x)
	\sqrt{g(x)}. \la{defA}
\end{align}
Then, we can write
\be
	E_C(u;\Omega) = \int_{\Omega} E^{\alpha\beta}_{ij}(x,u,Du)D_\alpha u^i D_\beta u^j dx.
	\la{EbyE}
\ee

In case that $m=\mbox{dim}(M)=2$,  the H\"{o}lder  continuity of a energy
minimizing map is shown in \cite{tac09}. 
For a energy minimizing map between Riemannian manifolds,
or more generally for a minimizer $u$ of a quadratic  functional
\[
	\int A^{\alpha \beta}_{ij}(x,u) D_\alpha u^i D_\beta u^j dx
\]
with smooth coefficients $A^{\alpha \beta}_{ij}(x,u)$,
once the H\"{o}lder continuity of $u$ has been shown, 
we see that
the coefficients $A^{\alpha \beta}_{ij}(x,u(x))$ are H\"{o}lder continuous,
and therefore we can show the $C^{1,\alpha}$-regularity of $u$
by virtue of Schauder-type estimate. Then, inductively we get 
higher regularity.
In contrast, if the target manifold is a Finsler manifold, 
the H\"{o}der continuity of $u$ does not imply the 
continuity of the coefficients $E^{\alpha\beta}_{ij}(x,u(x),Du(x))$.
So, if we want to obtain $C^{1,\alpha}$-regularity of a minimizer,
we have to show it directly.

In differential geometric setting, usually one assumes $C^\infty$-regualrity 
on the metric as (F-1). 
However, to get $C^{0,\alpha}$- or $C^{1,\alpha}$-regularity 
for energy minimizing maps,
it is enough to emply the following conditions instead of (F-1)
\begin{description}
	\item {(F-1a)} There exists a concave increasing function $\omega : [0, \infty)
	\to [0, \infty)$ with $\lim_{t\to +0} \omega(t)=0$ such that
	\begin{equation}
	 |F^2(u,X)-F^2(v,X)| \leq \omega(|u-v|^2) |X|^2 \la{omega-F}
	\end{equation}
	holds for any $u,v \in \R^n$ and $X\in \R^n$.
	\item {(F-1b)} $F(u,X)$ is twice differentiable in $X$ 
	for every $(u,X)\in T\R^n\setminus 0$.
\end{description}

On the other hand, about convexity we need the following 
uniformly convexity condition which is 
stronger than (F-3).
\begin{description}
	\item{(F-3a)}
	 There exist  positive constants $\lambda< \Lambda$
	for which
	\begin{equation}
     	  \lambda |\xi|^2 \leq f_{ij} (u,X)\xi^i\xi^j 
	 =\frac{1}{2} \frac{\partial^2 F^2(u,X)}{\partial X^i \partial X^j}
	 \xi^i \xi ^j \leq
     	\Lambda|\xi|^2
     	\la{ellip-f}
	\end{equation}
	holds for any $u,v \in \R^n$ and $(X,\xi) \in (\R^n \setminus 0) \times \R^n$.
\end{description}

The main result of this paper is as follows.
\begin{thm}\la{main}
Let $(M,g)$ a 2-dimentional smooth Riemannian manifold,
$\Omega \subset M$ a bounded domain with smooth boundary
$\partial \Omega$  and  $(\R^n,F)$ a Finsler space
with the Finsler structure $F$ satisfying
(F-1a), (F-1b), (F-2) and (F-3a).
Let $u\in H^{1,2}(\Omega, \R^n)$ be an energy minimizing map
in the class
\[
	H^{1,2}_\phi(\Omega, \R^n):=
	\{v \in H^{1,2} (\Omega, \R^n)~;~v-\phi 
	\in H^{1,2}_0(\Omega, \R^n)\}.
\]
Then $u\in C^{1,\alpha}(\Omega)
\cap C^{0,\beta}(\overline{\Omega})$ for some $\alpha\in (0,1)$
and any $\beta \in (0,1)$.
\end{thm}


\section{Proof of Theorem 1.1}
\setcounter{equation}{0}

In order to prove Theorem 1.1, we prepare the following higher integrability
results of minimizers which can be deduced easily from
\cite[Lemma 1]{josmei83} as mentioned in \cite{tac09}.
\begin{lem}[{\cite[Remark 5.3]{tac09}}]\la{Lq-est-Du}
Let $(M,g)$ be a smooth Riemannian m-manifold and
$\Omega \subset M$ a bounded domain with smooth boundary
$\partial \Omega$  and  $(\R^n,F)$ a Finsler space
with the Finsler structure $F$ satisfying
\pref{ellip-f}.
Suppose that $\phi \in H^{1,p}(\Omega, \R^n)$ for some $p>2$.
Let $u\in H^{1,2}(\Omega, \R^n)$ be an energy minimizing map
in the class $H^{1,2}_\phi(\Omega, \R^n)$.
Then, there exists a positive number $q_0>2$
such that for every $q\in (2,q_0)$, the estimate
\be\la{Lq-est-u}
	\int_\Omega |Du|^q dx \leq C \int_\Omega |D\phi|^q dx
\ee
holds.
\end{lem}

~~

Now, using several estimates which are obtained in \cite{tac09},
we can show the main result of this paper.
In  \cite{tac09}
the author supposed that 
\[
	A(x,u,p)=E^{\alpha\beta}_{ij}(x,u,p)p^i_\alpha p^j_\beta
\]
is in the class $C^{1,1}({\cal X}) \cap C^3({\cal X}^\prime)$, where
\[
	{\cal X}= \Omega \times \R^n \times \R^{mn} ~~\mbox{and}~~
	{\cal X}^\prime= \Omega \times \R^n \times (\R^{mn}\setminus \{0\}).
\]
However, it is clearly superfluous to obtain 
$C^{0,\alpha}$-regularity of the minimizer. 
In fact, it is easy to see that every proof in \cite{tac09} can be carried
assuming on the regularity of $A(x,u,p)$ only that
\begin{enumerate}
	\item $A(x,u,p)$ is in the class $C^{1,1}({\cal X})$ and
	 twice differentiable in $p$ at every
	$(x,u,p)\in {\cal X}^\prime$.
	\item There exists a concave increasing function
	$\omega:[0,\infty) \to [0,\infty)$ with $\lim_{t\to 0} \omega(t)=0$
	such that
	\[
	|A(x,u,p) - A(y,v,p)| \leq \omega (|x-y|^2+ |u-v|^2)|p|^2,~~
	\]
	holds for all
	$x, y \in \Omega, u,v \in \R^n$ and
	$p \in \R^{mn}\setminus 0 $.
\end{enumerate}
Therefore, all results in \cite{tac09} hold under the assumptions
in Theorem \ref{main} in the present paper.

If $u:\Omega \subset M \to \R^n$ minimizes the energy functional on $\Omega$,
then $u$ minimizes it on every sub-domain of $\Omega$. On the other hand,
the regularity  is a local property. So, it is suffices to study the regularity problem
on a domain $\Omega \subset \R^m$.

\begin{proof}[Proof of Therem \ref{main}]
First, we show that $u \in C^{0,\beta}(\overline{\Omega})$ for any $\beta \in (0,1)$.

We use the following notation as in \cite{tac09}.
For $x\in \Omega$ and $R>0$ we put
\begin{equation}\la{defQ}
	Q(x, R):= \{ y \in \R^m ~;~ |y^\alpha -x^\alpha| <R,~ \alpha=1, \ldots ,m\}.
\end{equation}	
For $x_0 \in \partial \Omega$ we always choose local coordinates so that
for sufficiently small $R_0>0$
\begin{align}
	Q(x_0,R_0) \cap \Omega \subset \R^m_+ = 
	\{x \in \R^m~;~ x^m>0\}, \nn\\
	Q(x_0,R_0) \cap\partial \Omega \subset 
	\{x \in \R^m~;~ x^m=0\}, \nn
\end{align}
and put for $0<R<R_0$
\begin{equation}\la{defQ+}
	Q^+(x_0, R):= Q(x_0, R) \cap \{x\in \R^m ~;~ x^m>0\}.
\end{equation}
Sometimes we write also 
\begin{equation}\la{defOmega+}
	\Omega(x,R):= \{y\in \Omega~;~|y^\alpha - x^\alpha| < R, ~\alpha =1,\ldots, m\},
\end{equation}
for general $x\in \Omega$ and $R>0$.

From \cite[(5.9)]{tac09}, when $x_0$ is an interior point and
$Q(x_0, 2r) \subset \subset \Omega$, we have for any $\delta \in (0,1)$ 
\begin{equation}\la{5.9}
	\begin{aligned}
		 & \int_{Q^(x_0,\rho)}|Du|^2 dx\\
		&  \leq  C \left\{ \left( \frac{\rho}{r} \right)^{2-\delta}
		 +\tilde{\omega} \big( r^2 +  \int_{Q(x_0,2r)} |Du|^2 dx \big)
		 \right\} \int_{Q^(x_0,2r)} |Du|^2 dx,
	\end{aligned}
\end{equation}
where $\tilde{\omega}=\omega^{(q-2)/q}$
for some $q>2$.
For a boundary point $x_0$, assuming
that $\phi \in H^{1,s} (s>m=2)$, from
\cite[(5.10)]{tac09}, we
have for any $\delta \in (0,1)$
\begin{equation}\la{5.10}
	\begin{aligned}
		 & \int_{Q^+(x_0,\rho)}|Du|^2 dx\\
		 & \leq C  \left\{ \left( \frac{\rho}{r} \right)^{2-\delta}
		 +\tilde{\omega} \big( r^2 + \int_{Q^+(x_0,2r)} |Du|^2 dx \big)
		 \right\} \int_{Q^+(x_0,2r)} |Du|^2 dx \\
		 & ~~~ + C(\phi) r^\gamma,
	\end{aligned}
\end{equation}
where $\gamma = 2(1-2/s)>0$.
Since we are assuming that $\phi \in H^{1,\infty}$, we can take 
$\gamma =2- \varepsilon$ for any $\varepsilon>0$.

Let us choose 
$\delta$ so that $2-\varepsilon < 2-\delta$.
Proceeding as in \cite[pp.317--318]{giu03},
we can deduce from \pref{5.9}
and \pref{5.10} that 
\begin{align}
	\int_{Q(x_0,\rho)}|Du|^2 dx &
	\leq M_1 \left( \frac{\rho}{r} \right)^{2-\varepsilon}
	 \int_{Q(x_0,r)} |Du|^2 dx ~~ \mbox{for} ~x_0 \in \Omega, \\
	\int_{Q^+(x_0,\rho)}|Du|^2 dx &
	\leq M_1 \left( \frac{\rho}{r} \right)^{2-\varepsilon} 
	 \int_{Q^+(x_0,r)} |Du|^2 dx +M_2 \rho^{2-\varepsilon} ~~\mbox{for}~ x \in \partial \Omega,
\end{align}	
for sufficiently small $r>0$ and $\rho \in (0,r)$,
where $M_1$ and $M_2$ are constants depending on
$g, F, \Omega$ and $\phi$.
Here, we used also the fact that
\begin{equation}\la{condforreg}
	\lim_{r_0 \to 0}
	\big\{ r_0^2+  \int_{\Omega(x_0,2r_0)}|Du|^2 dx \big\}
	=0
\end{equation}
holds for any $x_0 \in \overline{\Omega}$.

Now, proceeding as in \cite[pp.318--319]{giu03},
we can have that 
for any $\varepsilon \in (0,1)$ there exists a positive constant $M$ such that
\begin{equation}\la{Mor-est-Du}
		\int_{\Omega(x_0,\rho)}|Du|^2 dx \leq \rho^{2-\varepsilon} M,
\end{equation}
for any $x_0 \in \overline{\Omega}$.
So, putting $2\beta = 2-\varepsilon$, by Morrey's Dirichlet growth theorem,
we see that
$u \in C^{0,\beta}(\overline{\Omega})$.

Let us show $C^{1,\alpha}$-regularity of $u$, proceeding as in \cite{giagiu83}.
For a cube $Q_0=Q(x_0, R) \subset \subset \Omega$,
we consider the following frozen functional $A^0$ defined by
\be\la{frozen}
	A^0(v) =\int_{Q_0} E^{\alpha\beta}_{ij}(x_0, u_R, Dv)D_\alpha v^i D_\beta v^j dx,
\ee
where 
\[
	u_R= \intmean_{Q_0} u dx.
\]
Let $v$ be a minimizer of $A^0$ in the class
\[
	\{ v\in H^{1,2}(Q_0)~;~ v-u \in H^{1,2}_0(Q_0)\}.
\]
Since $u\in H^{1,q}$ for every $q\in (2,q_0)$ for some $q_0>2$ by 
Lemma \ref{Lq-est-Du}, using 
Lemma \ref{Lq-est-Du} for $v$, we see that there exists a positive number $q_1>2$
such that for every $q\in (2,q_1)$ there holds 
\be\la{Lq-est-Dv<Du}
	\int_{Q_0} |Dv|^q dx \leq \int_{Q_0} |Du|^q dx.
\ee
Moreover, as in \cite{tac09}, by using of difference quotient method, we can see that
$v \in H^{2,2}$ and that $Dv$ satisfies a system of 
uniformly elliptic equations weakly.
So,  for any $Q(x,r) \subset Q_0$, $Dv$ satisfies the Caccioppoli inequality,
\be\la{caccio}
	\int_{Q(x,r/2)} |D^2v|^2 dy \leq \frac{C}{r^2} \int_{Q(x,r)} |D-(Dv)_r|^2 dy,
\ee
and $D^2v$ satisfies \it reverse H\"{o}lder inequalities with increasing supports \rm 
due to Giaquinta-Modica (cf. \cite[p.299, Theorem 3]{giamodsou98-2},  
\be\la{rev-hol-v}
	\left( \intmean_{Q(x,r/2)} |D^2v|^q dy \right)^{1/q}
	\leq C \left( \intmean_{Q(x,r)} |D^2v|^2 dx\right)^{1/2},
\ee
for every $q\in (2,q_2)$ for some $q_2>2$.

Since we are considering 2-densional case, 
the Sobolev-Morrey imbedding theorem (cf. \cite[Theorem 3.11]{giu03}
yields that $v \in C^{1,\delta}$ for $\delta= 1-(2/q)$. Moreover, we have
for $\rho \in (0,R/4)$
\begin{align}
	& \left\{\rho^{-2-2\delta}\int_{Q(x_0,\rho)} |Dv-(Dv)_\rho|^2 dx \right\}^{1/2}\nn\\
	\leq & \sup_{Q(x_0,R/4)} \frac{|Dv(x)-Dv(y)|}{|x-y|^\delta}
	\leq C
	\|D^2v\|_{L^q(Q(x_0,R/4)}.\la {est-v-1}
\end{align}
For the last inequality,  we used Morrey-type inequality. 

Combining \pref{est-v-1}, \pref{rev-hol-v} and \pref{caccio}, we obtain
\begin{align}
	& \left\{\rho^{-2-2\delta}\int_{Q(x_0,\rho)} |Dv-(Dv)_\rho|^2 dx \right\}^{1/2}\nn\\
	\leq & C R^{\frac{2}{q}-1}\|D^2v\|_{L^2(Q(x_0, R/2)} \nn\\
	\leq & \left( R^{-2-2\delta} \int_{Q(x_0, R)} |Dv-(Dv)_R|^2 dx \right)^{1/2}.
	\la{est-v-2}
\end{align}

Putting $w=u-v$, we obtain
\begin{align}
	& \int_{Q(x_0,\rho)} |Du-(Du)_\rho|^2 dx \nn\\
	\leq & \int_{Q(x_0,\rho)} |Du-(Dv)_\rho|^2 dx \nn\\
	\leq & \int_{Q(x_0,\rho)} |Dv-(Dv)_\rho|^2 dx 
		+\int_{Q(x_0,\rho)} |Dw|^2 dx \nn\\
	\leq  & C\Big(\frac{\rho}{R}\Big)^{2+2\delta} \int_{Q(x_0, R)} |Dv-(Dv)_R|^2 dx 
		+C \int_{Q(x_0,\rho)} |Dw|^2 dx\nn\\
	\leq & C\Big(\frac{\rho}{R}\Big)^{2+2\delta} \int_{Q(x_0, R)} |Dv-(Du)_R|^2 dx 
		+C \int_{Q(x_0,\rho)} |Dw|^2 dx\nn\\
	\leq & C\Big(\frac{\rho}{R}\Big)^{2+2\delta} \int_{Q(x_0, R)} |Du-(Du)_R|^2 dx 
		+C \int_{Q(x_0,R)} |Dw|^2 dx.
		\la{est-Du-(Du)}
\end{align} 

Let us estimate $\int |Dw|^2dx$.
Proceeding as in \cite[pp.1967-1968]{tac09}, it is easey to see that
\begin{align}
	\int_{Q(x_0,R)} |Dw|^2 dx 
	\leq &  C \Big[ \int_{Q(x_0,R) } \omega(|x-x_0|^2 + |u-u_R|^2) |Du|^2 dx \nn\\
	&
	+  \int_{Q(x_0,R) }\omega(|x-x_0|^2 + |v-u_R|^2) |Dv|^2 dx\Big] \la{est-Dw-1}\\
	=:  & ~ I +II.\nn
\end{align} 
Using Jensen's inequality, H\"{o}lder's inequality and reverse H\"{o}lder ineqalitty, we
can estimate $I$ as follows.
\begin{align}
	I &  \leq C\Big( \int_{Q(x_0,R)} \omega^{q/(q-2)}dx \Big)^{(q-2)/q}
	\Big(\int_{Q(x_0,R)} |Du|^q \Big)^{2/q} \nn \\
	& \leq C \Big( \intmean_{Q(x_0,R)}\omega dx )^{(q-2)/q} 
	R^{m(q-2)/q}
	\Big( \int_{Q(x_0, R)} |Du|^q dx \Big)^{2/q} \nn\\
	& \leq C \Big( \omega \big( \intmean_{Q(x_0,R)} (|x-x_0|^2 + |u-u_R|^2) dx \big)
	\Big)^{(q-2)/q} R^{m(q-2)/q} R^{2m/q} \nn \\
	& ~~~~~~~~~~~~~ ~~~~~~~~~~~~~ ~~~~~~~~~~~ ~~~~~~~~~~~~~
	\cdot \Big(
	\intmean_{Q(x_0,R)} |Du|^q dx \Big)^\frac{2}{q} \nn\\
	& \leq C \Big( \omega \big( \intmean_{Q(x_0,R)} (R^2 + |u-u_R|^2) dx \big)
	\Big)^{(q-2)/q} \int_{Q(x_0,2R)} |Du|^2 dx. \la{est-Dw-I}
\end{align}
Here we used the boundedness of $\omega$.
By virtue of \pref{Lq-est-Dv<Du}, we can estimate $II$ similarly and get
\begin{align}
	II & \leq C\Big( \int_{Q(x_0,R)} \omega^{q/(q-2)}dx \Big)^{(q-2)/q}
	\Big(\int_{Q(x_0,R)} |Dv|^q \Big)^{2/q} \nn \\
	& \leq C \Big( \omega \big( \intmean_{Q(x_0,R)} (R^2 + |v-u_R|^2) dx \big)
	\Big)^{(q-2)/q} R^m
	\Big(\intmean_{Q(x_0,R)} |Du|^q dx \Big)^{2/q} \nn\\
	& \leq  C \Big( \omega \big( C \intmean_{Q(x_0,R)} 
	(R^2 + |u-u_R|^2+|v-u|^2) dx \big) 
	\Big)^\frac{q-2}{q}   \nn \\
	& ~~~~~~~~~~~~~~~~~~~~  ~~~~~~~~~ ~~~~~~~~~~~~~~~~~~~~\cdot
	 \int_{Q(x_0,2R)} |Du|^2 dx. \la{est-Dw-II}
\end{align}

Let us estimate the ingredients in $\omega$.
Using Sobolev's inequality (cf. \cite[p.103]{giu03}, 
we can see that for $2_\ast = 2m/(m+2)$
\begin{align}
	&  \intmean_{Q(x_0,R)} |u-u_R|^2 dx \nn \\
	\leq & C R^{-m} \Big( \int_{Q(x_0,R)} |Du|^{2_\ast} dx \Big)^{2/2_\ast}
	\nn \\
	\leq & C R^{-m} \Big( \int_{Q(x_0,R)} 1^{2/(2-2_\ast)} dx 
	\Big)^{2-2_\ast}
	\Big( \int_{Q(x_0,R)} |Du|^2 dx \Big) \nn\\
	\leq & C R^{-m + 2m -2_\ast m} \Big( \int_{Q(x_0,R)} |Du|^2 dx \Big) 
	\nn
\end{align} 
Since we are assuming that $m=2$, we  have $2_\ast =1$.
Thus, the above estimate together with \pref{Mor-est-Du} gives
for every $\varepsilon\in (0,1)$ the folowing estimate
\be\la{est-|u-uR|}
	\intmean_{Q(x_0,R)} |u-u_R|^2 dx
	\leq C \int_{Q(x_0,R)} |Du|^2 dx
	\leq C R^{2-\varepsilon}.
\ee
We can see also that
\begin{align}
	& \intmean_{Q(x\_,R)} |u-v|^2 dx \nn\\
	\leq & C \int_{Q(x_0,R)}\big(|Du|^2 + |Dv|^2\big) \nn \\
	\leq & C \int_{Q(x_0,R)} |Du|^2 dx \leq C R^{2-\varepsilon}
	\la{est-|u-v|}.
\end{align}
Since we can assume that $R \leq 1$, we see that
the ingredient in $\omega$ can be estimates by
$CR^{2-\varepsilon}$ for every $\varepsilon \in (0,1)$.

Using the assumption that $\omega (t) \leq C t^\sigma$ for 
some $\sigma \in (0,1]$, we obtain
\be\la{est-omega}
	\omega (...) \leq C R^{\sigma (2-\varepsilon)},
\ee
So, we can estimate $\omega^{(q-2)/q} \int |Du|^2 dx$
in \pref{est-Dw-I} and \pref{est-Dw-II} as 
\be\la{est-omega|Du|}
	\omega^{(q-2)/q} \int_{Q(x_0,2R)} |Du|^2 dx
	\leq CR^{(2-\varepsilon)\{1+\sigma(q-2)/q\} },
\ee 
where we used \pref{Mor-est-Du} again.
Now, take $\varepsilon\in (0,1)$ sufficiently small so that
\[
	(2-\varepsilon)\Big(1+\sigma \cdot \frac{q-2}{q} \Big)>2,
\]
and put
\be\la{est-omega|Du|-2}
	\gamma :=(2-\varepsilon)\Big(1+\sigma \cdot \frac{q-2}{q} \Big)-2>0.
\ee
Combining \pref{est-Dw-I}, \pref{est-Dw-II}, \pref{est-omega|Du|} and
\pref{est-omega|Du|-2}, we get
\be\la{est-Dw-4}
	\int_{Q(x_0, R)} |Dw|^2 dx
	\leq CR^{2+\gamma}.
\ee
Now, substituting the above inequality into
\pref{est-Du-(Du)}, we obtain
\begin{align}
	& \int_{Q(x_0,\rho)} |Du-(Du)_\rho|^2 dx \nn\\
	\leq  & C\Big(\frac{\rho}{R}\Big)^{2+2\delta} 
	\int_{Q(x_0, R)} |Du-(Du)_R|^2 dx 
		+ CR^{2+\gamma}. \la{est-Du-(Du)-2}.
\end{align}  

Using well known lemma (cf. \cite[Lemma 2.2]{giagiu83}, we conclude that
\be
	\int_{Q(x_0,\rho)} |Du-(Du)_\rho|^2 dx
	\leq C \rho^{2+2\alpha}
\ee
with $\alpha = \min \{\delta, \gamma/2\}$ for
every $Q(x_0, 2\rho) \subset \Omega$, and hence
$Du \in C^{\alpha}(\Omega)$.

\end{proof}

\begin{rem}
The perfect dominance functions treated by
S.Hildebrandt and H. von der Mosel
in \cite{hilvon03com, hilvon03jre}
have the structure similar to that of the energy density $e_c$. 
So, some of their results are valid for weakly harmonic maps in
2-dimensional case.
More precisely, for the case that $F(u,X)$ is 
continuously differentiable in $u$, 
once the H\"{o}der continuity of a weakly harmonic map have shown, 
we can get  its $C^{1,\alpha}$-regularity proceeding exactly
as in the fourth section of \cite{hilvon03com}.
On the other hand, in this paper,
we prove $C^{1,\alpha}$-regularity using the minimality
without assuming
the differentiability of 
$F(u,X)$ with respect to $u$.

We should mention also that
in \cite{hilvon03com} the minimality is not necessary 
to get $C^{1,\alpha}$-regularity for H\"{o}lder continuous
weak solutions of the Euler-Lagrange equation
of a perfect dominance function. 
However, in both of \cite{hilvon03com} and this paper,
the minimality is necessary to get the H\"{o}lder continuity.
\end{rem}

\end{document}